\documentclass[journal ]{new-aiaa}
\usepackage[utf8]{inputenc}
\usepackage{textcomp}

\usepackage{graphicx}
\usepackage{amsmath}
\usepackage[version=4]{mhchem}
\usepackage{longtable,tabularx}
\usepackage{algpseudocode}
\usepackage{algorithm}

\usepackage{physics}
\usepackage{multirow}
\usepackage{subfig}
\usepackage{comment}
\usepackage[title,titletoc]{appendix}

\usepackage{booktabs}
\allowdisplaybreaks

\title{Orbital Depot Location Optimization for Satellite Constellation Servicing with Low-Thrust Transfers}

\author{Euihyeon Choi \footnote{Ph.D. Student, Daniel Guggenheim School of Aerospace Engineering.} and Koki Ho \footnote{Associate Professor, Daniel Guggenheim School of Aerospace Engineering; kokiho@gatech.edu. Associate Fellow AIAA. (Corresponding Author)}}
\affil{Georgia Institute of Technology, Atlanta, Georgia 30332}

\begin{document}

\maketitle

\begin{abstract}
This paper addresses the critical problem of co-optimizing the optimal locations for orbital depots and the sequence of in-space servicing for a satellite constellation. While most traditional studies used network optimization for this problem, assuming a fixed set of discretized nodes in the network (i.e., a limited number of depot location candidates), this work is unique in that it develops a method to optimize the depot location in continuous space. The problem is formulated as mixed-integer nonlinear programming, and we propose a solution methodology that iteratively solves two decoupled problems: one using mixed-integer linear programming and the other using nonlinear programming with an analytic transfer solution. To demonstrate the effectiveness of our approach, we apply this methodology to a case study involving a GPS satellite constellation. Numerical experiments confirm the stability of our proposed solutions.
\end{abstract}

\section*{Nomenclature}

% \noindent(Nomenclature entries should have the units identified)

{\renewcommand\arraystretch{1.0}
\noindent\begin{longtable*}{@{}l @{\quad=\quad} l@{}}
$D_0$                   & index set of originally given depots \\
$D_e$                   & index set of end depots \\
$D_s$                   & index set of start depots \\
$D_v$                   & index set of virtual depots \\
$D_s^k$                 & index set of original and virtual depots corresponding to depot $k$ \\
$I_{sp}^d$              & specific impulse of the depot \\
$I_{sp}^l$              & specific impulse of the launch vehicle \\
$I_{sp}^s$              & specific impulse of the servicer \\
$i$                     & inclination of the orbit \\
$m^{\text{dry},d}_k$    & dry mass of depot $k$ \\
$m^{\text{dry},s}_k$    & dry mass of servicer $k$ \\
$m^L_i$                 & payload mass for satellite $i$ \\
$m^l_{\text{max}}$      & maximum launch mass per depot \\
$n_d$                   & number of depots \\
$n_t$                   & number of satellites \\
$n_v$                   & number of possible routes per depot \\
$r_0$                   & virtual radius for EMLEO calculation \\
$r_{\text{min}}$        & minimum possible radius of depot \\
$T$                     & index set of satellites \\
$u_k$                   & continuous decision variable indicating mass of vehicle $k$ after payload deployment \\
$z_{ij}^k$              & binary decision variable (1 if edge from node $i$ to $j$ is used by servicer $k$ and 0 otherwise) \\
$\Delta V^d$            & incremental velocity of the depot \\
$\Delta V^l$            & incremental velocity of the launch vehicle \\
$\Delta \theta$         & angle difference between two orbital planes \\
$\vb*{\eta}_k$          & orbital element vector of depot $k$ \\
$\mu$                   & standard gravitational parameter of the Earth \\
$\Omega$                & right ascension of the ascending node of the orbit \\
$\phi_k$                & EMLEO conversion factor for depot $k$ \\

\multicolumn{2}{@{}l}{Subscripts}\\
ECI     & Earth-centered inertial frame\\
$i$/$j$ & index representing a satellite \\
$k$     & index representing a depot or corresponding vehicle \\
peri    & perifocal frame\\

\multicolumn{2}{@{}l}{Superscripts}\\
$d$     & value related to the depot\\
$k$     & index representing the depot or corresponding vehicle \\
$L$     & value related to the payload\\
$l$     & value related to the launch vehicle\\
$s$     & value related to the servicer\\
\end{longtable*}}

\section{Introduction}
\lettrine{T}{he} increasing complexity of interconnected space systems demands innovative optimization frameworks to manage in-space operations effectively. As the industry shifts from isolated, single-satellite missions to constellations of collaborating satellites \cite{Marcuccio_2019_SmallSatTrend, Kulu_2024_ConstellationSurvey}, in-space servicing, assembly, and manufacturing (ISAM) technologies have become essential for reducing costs and extending satellite lifespans \cite{Saleh_2002_OOS1, Lamassource_2002_OOS2, Long_2007_OOS}. This paradigm shift not only drives the need for advanced operational strategies but also motivates the development of methodologies that optimize both initial deployments and ongoing maintenance of space assets. Consequently, the field of space logistics focused on designing and planning resource flows for long-term missions is playing a crucial role in enabling future space endeavors \cite{Ho_2024_SpaceLogisticsReview}.

The challenge of optimizing in-space servicing missions for multiple satellites can be effectively modeled as an extension of the well-known optimal routing problem. Traditionally, the problem of visiting all given multiple tasks (satellites) with multiple routes while minimizing total cost is known as the vehicle routing problem (VRP) \cite{VehicleRouting_book_ch2}. Numerous studies further extend the VRP with various objectives and constraints, such as collaborating multiple vehicles \cite{Choi_2024_interagent, Choi_2025_VGRPP}, considering synergy effects between exploration sites \cite{Choi_2022_VRPPCS}, and utilizing both rover and helicopter on the Martian surface \cite{Choi_2023_2EVRPP}. For Earth-orbiting systems application, Bang and Ahn \cite{Bang_2019_MultitargetRendezvous} have also studied the optimal routing problem in Earth-orbiting scenarios for active debris removal missions, referred to as multitarget rendezvous. Lee and Ahn \cite{Lee_2023_MultitargetRendezvousHybrid} have further extended the problem to utilize a hybrid propulsion system for transfer between orbits.

An effective ISAM strategy requires not only optimizing servicing sequences but also strategically determining facility (depot) locations to maximize mission efficiency. The facility location problem (FLP) determines the optimal facility (depot) locations to minimize the sum of costs while satisfying demand constraints of destinations \cite{Cooper_1963_FLP}. 
In space applications, Zhu et al. \cite{Zhu_2020_fuelStationOD} proposed FLP for asteroid mining and ISAM depots in a sun-synchronous orbit. Recently, Shimane et al. \cite{shimane_orbital_2024} presented an orbital FLP for satellite constellation servicing in high-altitude orbit with low-thrust propulsion systems. 

Combining depot placement with route optimization, the location routing problem (LRP) offers a powerful tool for addressing ISAM challenges. The LRP, a synthesis problem of FLP and VRP, typically assumes discrete depot locations and only a few studies have been conducted that allow continuous locations \cite{Salhi_2009_planarLRP, Manzour_2012_PlanarLRP, Mohammad_2021_planarLRP}. The LRP determines which depots to use from a set of candidates and simultaneously plans the routes to minimize overall costs \cite{Prodhon_2014_LRPsurvey, Mara_2021_LRPreview}. For space exploration context, Ahn et al. \cite{Ahn_2012_GLRPP} introduced a generalized LRP with profits to optimize depot strategies, routing tactics, and site selection using a column generation heuristic.

While the aforementioned studies have examined the depot location(s) and servicer routing problem, a common limitation is that they typically considered a fixed network to perform the optimization. This means that they need to start with a set of fixed nodes used as candidate depot locations or staging points and optimize the depot location(s) from those candidates. Since the size of the optimization problem is directly related to the number of nodes, these existing methods could only consider a relatively small number of discrete candidates for the depot location(s). Thus, the solution from these methods could be suboptimal. There is a need to optimize the depot location(s) in a continuous space for more rigorous and effective planning.

In response to this need, this paper introduces a continuous location routing problem for orbital depots (CLRPOD), extending the continuous LRP to the Earth-orbiting systems. Our approach can identify both depots' orbital elements and the visiting routes of each servicer, while minimizing total propellant consumption with practical constraints, such as maximum launch mass. A key difficulty arises because the transfer cost between two orbits cannot be predetermined, as the depots’ orbital elements are variables in the optimization. Moreover, since payload deliveries affect the servicer spacecraft's mass and thus simply minimizing the sum of \(\Delta V\) does not minimize the total propellant consumption. To address these challenges, our methodology decouples depot location optimization from route planning. We begin by fixing initial depot locations to determine optimal visiting sequences; then, using the routing solution, we iteratively calibrate the depots’ orbital elements until convergence. Note that we assume the low-thrust transfer between orbits, Edelbaum's analytic solution \cite{Edelbaum_1961_LowThrust, Edelbaum_1962_TheoryMaxMin} can accelerate the depots' orbital element optimization phase without numerical error.

There are three key elements in this paper.
First, we present a novel problem formulation and solution methodology that simultaneously optimizes orbital elements for multiple depots and servicing routes, with depot locations determined over a continuous space rather than from predetermined discrete candidates.
Second, the rocket equations are directly integrated into the combinatorial optimization problem to reflect realistic in-space servicing scenarios, such as delivering physical payloads to multiple satellites.
Lastly, we employ analytic derivatives of total propellant mass using low-thrust transfer, enabling fast and accurate optimization of depot orbital elements.

The structure of the paper is organized as follows. Section \ref{Sec:Problem_description} provides the definition and mathematical formulation of the proposed problem, including an analytical solution for low-thrust orbital evolution and a cost estimation for depots. Section \ref{Sec:Solution_methodology} details the methodology used to solve the complex mixed-integer nonlinear problem (MINLP). In Section \ref{Sec:Case_study}, we present a case study involving a constellation of GPS satellites along with numerical experimental results. Finally, Section \ref{Sec:Conclusions} summarizes the conclusions drawn from this work.

\section{Problem Description} \label{Sec:Problem_description}

\subsection{Continuous Location-Routing Problem for Orbital Depots}

The CLRPOD is a special application of LRP that refers to the problem of simultaneously optimizing the routes for the in-space servicing of a constellation of satellites and the orbital elements of orbital depots. The problem aims to utilize multiple depots to visit all given servicing demand satellites while minimizing the fuel required from an effective mass to LEO (EMLEO) perspective. The EMLEO is a hypothetical mass including the propellant required to transfer the depot from a (virtual) circular LEO to its desired orbit for standardized metric \cite{shimane_orbital_2024}. The servicers depart from each depot, visit multiple satellites, and return to the depot from which they departed. The problem allows multiple routes for each servicer, which can also be interpreted as deploying multiple servicers per depot. It is assumed that the servicers use low thrusters for orbit transitions, and the amount of fuel required for this can be calculated based on the derivation in Sec.~\ref{Sec:Low-Thrust Transfer Between Circular and Tilted Circular Orbits}. Note that the fuel mass should be converted to EMLEO by multiplying the conversion factor derived in Sec.~\ref{Sec:Depot Launch and Deployment Costs}, and the sum of these EMLEOs is minimized by considering both the depot's orbital elements and servicing routes. For simplicity, all demands are assumed to be known in advance for a static problem. Notably, the orbital elements of the depots are not discretized but are allowed to take continuous values within specified ranges. Furthermore, the highly nonlinear and complex impact of orbital element variations on the solution necessitates formulating this problem as a MINLP problem, which is inherently challenging to solve comprehensively.

It is important to note that this study focuses exclusively on regular, non-urgent maintenance and refueling operations, where service schedules can be fully planned in advance. In other words, mission time is not a critical constraint in the proposed problem. This assumption enables low-thrust propulsion systems to reduce overall fuel consumption, which aligns with minimizing EMLEO. In contrast, time-critical servicing missions such as emergency response or dynamic tasks require time constraints in the optimization and may involve high-thrust propulsion or other rapid deployment methods.

\begin{figure}[hbt!]
    \centering
    \subfloat[]{
        \includegraphics[width=0.45\textwidth]{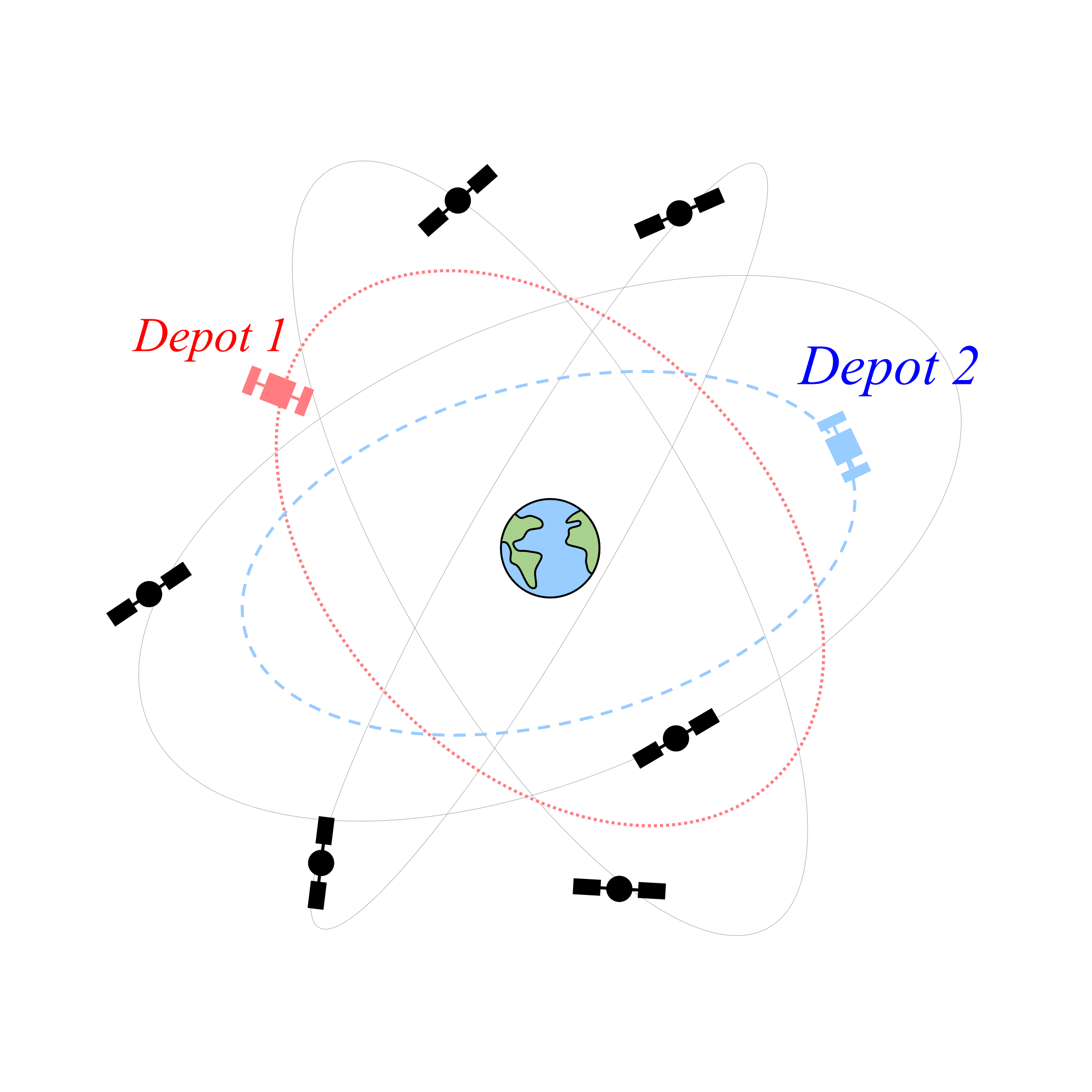}
        \label{fig:prob_def1}
    }
    \hfill
    \subfloat[]{
        \includegraphics[width=0.45\textwidth]{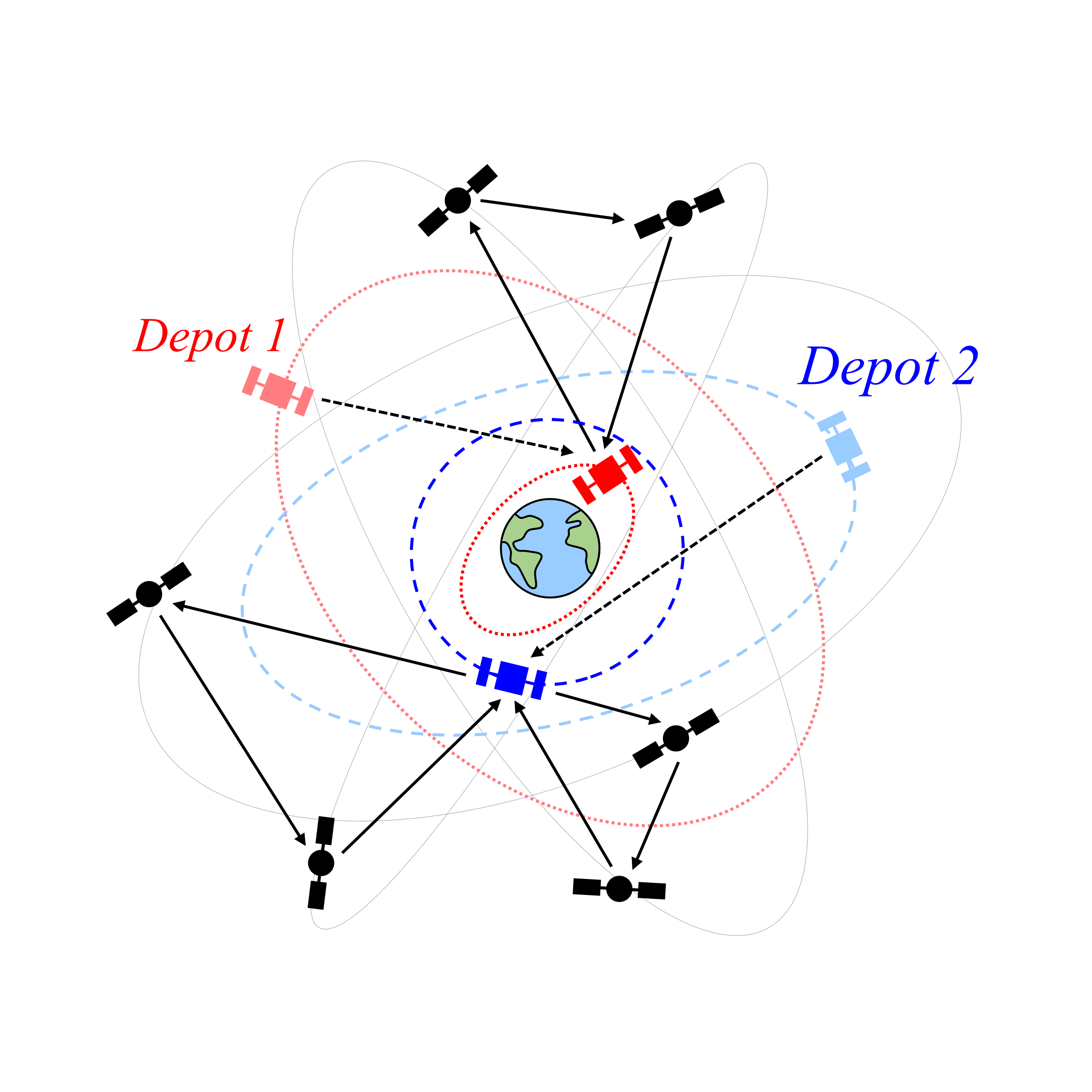}
        \label{fig:prob_def2}
    }
    \caption{A typical example of the CLRPOD: a) problem instance with initial orbits of orbital depots; b) routing solution with modified orbits of orbital depots.}
    \label{fig:prob_def}
\end{figure}

Figure \ref{fig:prob_def} illustrates a typical example of a solution for the CLRPOD. The scenario involves six satellites serviced by servicer spacecrafts from two depots on the initial orbits (Fig.~\ref{fig:prob_def1}). Figure \ref{fig:prob_def2} shows the optimal solution, including the servicing sequences (or routes) and the final orbits of the depots. At \textit{Depot 1}, the servicer operates along a single route. In contrast, at \textit{Depot 2}, the servicer initially visits two satellites, returns to \textit{Depot 2} to resupply, and then services the remaining two satellites. For illustrative purposes, transfers between orbits are depicted as simple arrows, while actual trajectories follow optimized low-thrust transfer paths. The orbital elements of \textit{Depots 1 and 2} (i.e., semimajor axis, inclination, right ascension of the ascending node, eccentricity, and argument of perigee) are also optimized for minimizing the propellant EMLEOs. For simplicity, this study assumes circular orbits, excluding eccentricity and argument of perigee from the optimization process. Nonetheless, the methodology presented in this paper can be extended to general elliptical orbits.

\subsection{Low-Thrust Transfer Between Circular and Tilted Circular Orbits} \label{Sec:Low-Thrust Transfer Between Circular and Tilted Circular Orbits}
Since the demands of satellite constellations are typically known in advance, low-thrust transfers are widely employed in on-orbit servicing to save fuel, extend servicer lifetimes, and reduce overall mission costs. One characteristic of low-thrust transfers is that the cost of phase maneuvers is relatively negligible due to the extended transfer times, with most transfer costs arising from changes in the orbit's shape and orientation in 3D.

Accurately calculating propellant expenditure for low-thrust transfers requires solving a complex trajectory optimization problem given specific initial and final boundary conditions, a process that can be computationally intensive. Consequently, a Lyapunov-function-based feedback control law known as the Q-law is widely used despite its suboptimal performance because it significantly reduces computation time. However, even the Q-law can be impractical for high-level decision-making when numerous iterations are needed.

The design of satellite constellations--characterized by similar orbital elements such as semimajor axis, eccentricities, and argument of perigee--not only ensures nearly uniform perturbation effects, thereby preserving orbital geometry and reducing station-keeping requirements, but also enables the application of analytical optimal solutions for circle-to-inclined-circle low-thrust transfers. This uniformity, especially in circular orbits where a constant altitude maintains consistent signal strength, is exemplified by formations like the Walker Delta Pattern constellation. Therefore, we can apply an analytical optimal solution for circle-to-inclined-circle low-thrust transfers. This solution, originally derived by Edelbaum \cite{Edelbaum_1961_LowThrust, Edelbaum_1962_TheoryMaxMin} and later reformulated by Kechichian \cite{Kechichian_1997_Reformulation}, assumes constant thrust acceleration, low thrust, and a quasi-circular transfer. 

Consider a servicer transferring from orbit 1 to orbit 2. The tilt angle between the two orbits can be defined as the angle between the normal vectors of the respective orbital planes. In other words, this angle can be measured by the angle between the unit vectors along the z-axes of each orbit’s perifocal frame, transformed into the Earth-Centered Inertial (ECI) frame (i.e., angle between $\hat{\vb*{z}}_{1, \text{ECI}}$ and $\hat{\vb*{z}}_{2, \text{ECI}}$, respectively).
Assuming circular orbits, $\hat{\vb*{z}}_{1, \text{ECI}}$ and $\hat{\vb*{z}}_{2, \text{ECI}}$ can be represented as follows.

\begin{align}
    \hat{\vb*{z}}_{1, \text{ECI}} &= R_z(-\Omega_1) R_x (-i_1) 
    \begin{bmatrix} 0 & 0 & 1 \end{bmatrix}_{\text{peri}_1}^\top = 
    \begin{bmatrix} \sin \Omega_1 \sin i_1 & -\cos \Omega_1 \sin i_1 & \cos i_1 \end{bmatrix}_{\text{ECI}}^\top \\
    \hat{\vb*{z}}_{2, \text{ECI}} &= R_z(-\Omega_2) R_x (-i_2) 
    \begin{bmatrix} 0 & 0 & 1 \end{bmatrix}_{\text{peri}_2}^\top = 
    \begin{bmatrix} \sin \Omega_2 \sin i_2 & -\cos \Omega_2 \sin i_2 & \cos i_2 \end{bmatrix}_{\text{ECI}}^\top
\end{align}
where $R_z(-\Omega_1)$/$R_z(-\Omega_2)$, as well as $R_x (-i_1)$/$R_x (-i_2)$ are the rotation matrices about the z-axis and x-axis with angle $-\Omega_1$/$-\Omega_2$ and $-i_1$/$-i_2$, respectively. Note that $\Omega_1$ and $\Omega_2$ represent the right ascension ascending node (RAAN), and $i_1$ and $i_2$ denote the inclination of orbit 1 and orbit 2, respectively.

Therefore, the tilt angle between the two orbits, in terms of their orbital elements, is given by

\begin{align}
    \Delta \theta &= \cos^{-1} \left( \hat{\vb*{z}}_{1, \text{ECI}} \cdot \hat{\vb*{z}}_{2, \text{ECI}} \right) \\
    &= \cos^{-1} \left( \sin \Omega_1 \sin i_1 \sin \Omega_2 \sin i_2 + \cos \Omega_1 \sin i_1 \cos \Omega_2 \sin i_2 + \cos i_1 \cos i_2 \right)\\
            &= \cos^{-1} (\sin i_1 \sin i_2 \cos (\Omega_1 - \Omega_2) + \cos i_1 \cos i_2) \label{Eq:dtheta}
\end{align}

The total $\Delta V$ required for the orbit transfer with low thrust is given by

\begin{equation}
    \Delta V = \sqrt{V_1^2 - 2 V_1 V_2 \cos \qty[\frac{\pi}{2} \min \qty{ \Delta \theta, 2} ] + V_2^2}
    \label{Eq:Edelbaum_dV_min}
\end{equation}
where $V_1$ and $V_2$ are the circular orbit velocities for orbits 1 and 2, respectively, defined as $\sqrt{\mu / r_1}$ and $\sqrt{\mu / r_2}$, where $r_1$ and $r_2$ are the radii of circular orbits 1 and 2. 
The full derivation of Eq.~\eqref{Eq:Edelbaum_dV_min} is detailed in previous literature \cite{Edelbaum_1961_LowThrust, Kechichian_1997_Reformulation}. The unit of the inclination change, $\Delta \theta$, is radians. The factor of $\frac{\pi}{2}$ appears because it connects the desired tilt angle (relative inclination) change ($\Delta \theta$) with the yaw steering of the spacecraft, reflecting the average effectiveness of the out-of-plane thrust component when applied continuously throughout the orbital plane change. Note that if $\Delta \theta > 2$, then $\Delta V = V_1 + V_2$ \cite{Kechichian_1997_Reformulation}.

\subsection{Depot Launch and Deployment Costs} \label{Sec:Depot Launch and Deployment Costs}
Accurate estimation of launch and insertion costs for orbital depots requires the consideration of both their masses and orbital elements, not merely summing depot masses. For example, if two depots have the same mass but are positioned at different altitudes, the one at a higher altitude incurs a greater launch cost. To address this challenge, we use the concept of EMLEO, which assumes each depot is transferred from a virtual circular orbit with radius $r_0$ to its intended orbit via a coplanar Hohmann transfer. The EMLEO is determined by adding the propellant mass required for the transfer. Additionally, we assume that the initial (hypothetical) burn is performed by the launch vehicle while the orbital depot executes the second burn. Mathematically, let $r_0$ be the radius used in calculating EMLEO and $a_k$ be the orbital depot $k$'s radius. Then, the required $\Delta V$ by the launch vehicle ($\Delta V^l$) and depot's deployment ($\Delta V^d$) are calculated as follows.
\begin{align}
    \Delta V^l(a_k) = \sqrt{\frac{2 \mu}{r_0} - \frac{2 \mu}{a_k + r_0}} - \sqrt{\frac{\mu}{r_0}} \\
    \Delta V^d(a_k) = \sqrt{\frac{\mu}{a_k}} - \sqrt{\frac{2 \mu}{a_k} - \frac{2 \mu}{a_k + r_0}}
\end{align}
where $\mu$ is the gravitational parameter of the Earth. Then, the EMLEO conversion factor for depot $k$ ($\phi_k$) is defined as

\begin{align}
    \phi_k(a_k) = \exp \left[\frac{\Delta V^l(a_k)}{g_0 I_{sp}^l}\right] \exp \left[ \frac{\Delta V^d(a_k)}{g_0 I_{sp}^d}\right]
    = \exp \left[ 
    \frac{ \sqrt{\frac{\mu}{a_k}} - \sqrt{\frac{2 \mu}{a_k} - \frac{2 \mu}{a_k + r_0}}}{g_0 \ I_{sp}^d} + \frac{\sqrt{\frac{2 \mu}{r_0} - \frac{2 \mu}{a_k + r_0}} -\sqrt{\frac{\mu}{r_0}}}{g_0 I_{sp}^l} \right] \label{eq:phi_k(a_k)}
\end{align}
where $g_0$ is the standard acceleration of gravity (i.e., $g_0 \approx 9.81$ m/s\textsuperscript{2}), $I_{sp}^l$ and $I_{sp}^d$ are the specific impulse of the launch vehicle and the depot, respectively. 

Note that this paper does not model the launch cost from Earth's surface to the parking orbit. The primary focus is optimizing depot placement and the servicing sequence among known satellites. Incorporating detailed launch vehicle trajectory design and deployment options (e.g., dedicated, rideshare, piggyback missions, launch vehicle selection, or launch site choice) is beyond the current scope and could be explored in future work. Moreover, the EMLEO metric includes the transportation cost from parking orbit to the final servicing orbit, which allows for a meaningful optimization of in-space operations and propellant usage.

\subsection{Mathematical Formulation of Orbital Depot Location-Routing Problem}
Let the index set of original depots ($D_0$), virtual depots ($D_v$), start nodes ($D_s$), end nodes ($D_e$), and satellites ($T$) be as follows.
\begin{align}
    D_0 &= \{0, 1, ..., n_d - 1\} \\
    D_v &= \{n_d, ..., n_d n_v - 1\} \\
    D_s &= D_0 \cup D_v \\
    D_e &= \{n_d n_v, ..., n_d (n_v + 1) - 1 \} \\
    T &= \{n_d (n_v + 1), ..., n_d (n_v + 1) + n_t - 1 \}
\end{align}
where $n_d$ is the number of depots, $n_v$ is the number of available vehicles (routes) for each depot, and $n_t$ is the number of satellites. Note that the vehicles associated with each depot can be modeled as multiple depots, including the original depot and its duplications (virtual depots), for mathematical convenience. For example, if $n_v = 1$, only original depots exist, meaning each depot has only a single route available. 

To manage the formulation involving virtual depots, we define two auxiliary indices for each start node $k \in D_s$. First, $k'(k) \in D_0$ denotes the original depot associated with the start node $k$, computed as $k'(k) \equiv \text{mod}(k, n_d)$. Second, $k''(k) \in D_e$ represents the corresponding end node, defined as $k''(k) \equiv n_d n_v + \text{mod}(k, n_d)$. Both $k'(k)$ and $k''(k)$ are thus functions of the start node index $k \in D_s$.
The set of original and virtual depots corresponding to an original depot $k' \in D_0$ is defined as follows:

\begin{align}
    D_s^{k'} = \{k', k'+n_d, ..., k'+n_d(n_v-1)\} \qquad \forall k' \in D_0
\end{align}

Now, we propose the MINLP for the satellite constellation servicing as follows:

($P_0$) Original CLRPOD formulation
\begin{align}
    \min_{\vb*{z}, \vb*{u}, \vb*{\eta}} \left[ \sum_{k \in D_s} \left( u_k - m_{k'(k)}^{\text{dry},s} \sum_{j \in T} z_{kj}^k \right) \phi_{k'(k)}(\vb*{\eta}_{k'(k)})\right] \label{eq:obj}
\end{align}
subject to

\begin{align}
    &\sum_{j \in T} z_{kj}^k \leq 1 \qquad \forall k \in D_s \label{eq:const1} \\ 
    &\sum_{k \in D_s} \left( z_{kj}^k + \sum_{i \in T, i \neq j} z_{ij}^k \right) = 1 \qquad \forall j \in T \label{eq:const2} \\
    &z_{kj}^k - z_{jk''(k)}^k + \sum_{i \in T, i \neq j} \left( z_{ij}^k - z_{ji}^k \right) = 0 \qquad \forall k \in D_s, j \in T \label{eq:const3}\\
    &u_k + M \left( 1 - z_{kj}^k \right) \geq \left( u_j + m^{L}_j \right) \exp \left( \frac{\Delta V_{k'(k) j}(\vb*{\eta}_{k'(k)})}{g_0 I_{sp}^s} \right) \qquad \forall k \in D_s, j \in T \label{eq:const4}\\
    &u_i + M \left( 1 - \sum_{k \in D_s} z_{ij}^k \right) \geq \left( u_j + m^{L}_j \right) \exp \left( \frac{\Delta V_{ij}}{g_0 I_{sp}^s} \right) \qquad \forall i \in T, j \in T, i \neq j \label{eq:const5}\\
    &u_i + M \left( 1 - z_{ik''(k)}^k \right) \geq m^{\text{dry},s}_{k'(k)} \exp \left( \frac{\Delta V_{i k'(k)}(\vb*{\eta}_{k'(k)})}{g_0 I_{sp}^s} \right) \qquad \forall k \in D_s, i \in T \label{eq:const6}\\
    &\left[\sum_{p \in D_s^{k'}} \left( u_p - m_{k'}^{\text{dry},s} \sum_{j \in T} z_{pj}^p \right) + m_{k'}^{\text{dry},s} + m_{k'}^{\text{dry},d} \right] \phi_{k'}(\vb*{\eta}_{k'}) \leq m_{\text{max}}^{l} \qquad \forall k' \in D_0\label{eq:const7}\\
    &z_{kj}^k \in \{0, 1\} \qquad \forall k \in D_s, j \in T \label{eq:const8}\\
    &z_{ij}^k \in \{0, 1\} \qquad \forall k \in D_s, i \in T, j \in T, i \neq j \label{eq:const9}\\
    &z_{ik''(k)}^k \in \{0, 1\} \qquad \forall k \in D_s, i \in T \label{eq:const10}\\
    &u_i \geq 0 \qquad \forall i \in D_s \cup T \label{eq:const11}
 \end{align}

Equation \eqref{eq:obj} is the objective function, which minimizes the overall propellant mass in terms of EMLEO. $u_k$ is the mass before departure from depot $k$ (i.e., sum of the servicer's dry mass, the required propellant mass, and the mass of all payloads to be serviced), $m_{k'}^{\text{dry},s}$ is the dry mass of servicer $k'$, $z_{ij}^k$ is the binary decision variable, which equals one if the vehicle $k$ travels from node $i$ to $j$ and zero otherwise, and $\phi_{k'}$ is the EMLEO conversion factor to convert the propellant mass into EMLEO. $\vb*{\eta}_{k'} = [a_{k'}, i_{k'}, \Omega_{k'}]$ is the set of orbital elements of depot $k'$. Since $\phi_{k'}$ is a nonlinear function of the semimajor axis of depot $k'$ ($a_{k'}$), the objective function is nonlinear.

The use of the EMLEO conversion factor ensures that the objective function inherently accounts for the delivery cost of both depot deployment and servicing operations. The physical meaning of the objective is the effective propellant mass required to conduct all servicing missions and deliver them in their respective orbits. Importantly, the model does not assume that 1 kg of propellant used by a servicer is equivalent to 1 kg used by a launch vehicle or depot. Instead, each mission stage (i.e., depot delivery and satellite servicing) is modeled using its respective propulsion parameters and $\Delta V$ requirements. Consequently, the EMLEO conversion factor appropriately combines launch, deployment, and servicing fuel consumption into a single objective value.

Equation \eqref{eq:const1} indicates that every start node (original or virtual depot) should be used at most once. 
Equation \eqref{eq:const2} enforces every given satellite should be visited exactly once, and Eq.~\eqref{eq:const3} represents the flow continuity constraint. 
Equations (\ref{eq:const4}--\ref{eq:const6}) ensure that the mass at each node reflects the transfer propellant mass calculated by the rocket equations. $u_i$ is the mass of servicer at node $i$ after deploy the payload ($m_i^L$), $M$ is a sufficiently large positive constant (big-M), $\Delta V_{ij}$ is the required velocity increment for orbit transfer from node $i$ to node $j$, and $I_{sp}^s$ is the specific impulse of the servicer.
Equation \eqref{eq:const7} restricts the maximum launch mass for each depot in terms of the EMLEO. Note that $\Delta V_{k' j}$, $\Delta V_{i k'}$, and $\phi_{k'}$ are the functions of $\vb*{\eta}_{k'}$, which make constraints nonlinear. 
Equations (\ref{eq:const8}--\ref{eq:const10}) and Eq.~\eqref{eq:const11} impose the binary and nonnegative constraints for the decision variables.
For notational convenience, we will write $k'(k)$ as $k'$ and $k''(k)$ as $k''$ for the rest of this paper.

\section{Solution Methodology} \label{Sec:Solution_methodology}
To effectively solve the CLRPOD defined in the previous section, this section introduces a methodology that alternates between solving two decomposed subproblems: mixed-integer linear programming (MILP) and nonlinear programming (NLP). Solving these subproblems alternately to address a complex MINLP is a common approach in the literature \cite{isaji_multidisciplinary_2022}, as each subproblem is simpler and has specialized solvers compared to the original MINLP. When the orbital elements of each orbital depot are fixed, the CLRPOD reduces to an MILP. Similarly, when the route of each servicer is fixed, the CLRPOD reduces to the NLP. By alternately solving these subproblems, we obtain a suboptimal solution to the CLRPOD.

\subsection{Mixed-Integer Linear Programming with Fixed Depots' Locations}
Assume all orbital elements of depots are fixed. Then, $\phi_{k'}$, $\Delta V_{k' j} (\vb*{\eta}_k')$, and $\Delta V_{i k'} (\vb*{\eta}_k')$ are all constants. So, the problem becomes an extension of the VRP, which can be formulated as follows:

($P_L$) Linearized $P_0$ with fixed depots' locations
\begin{align}
    \min_{\vb*{z}, \vb*{u}} \left[ \sum_{k \in D_s} \left( u_k - m_{k'}^{\text{dry},s} \sum_{j \in T} z_{kj}^k \right) \phi_{k'}\right] \label{eq:linear obj}
\end{align}
subject to

\begin{center}
    Equations (\ref{eq:const1}--\ref{eq:const3}), \eqref{eq:const5}, and (\ref{eq:const8}--\ref{eq:const11})
\end{center}
\begin{align}
    &u_k + M \left( 1 - z_{kj}^k \right) \geq \left( u_j + m^{L}_j \right) \exp \left( \frac{\Delta V_{k' j}}{g_0 I_{sp}^s} \right) \qquad \forall k \in D_s, j \in T \label{eq:linear const4}\\
    &u_i + M \left( 1 - z_{ik''}^k \right) \geq m^{\text{dry},s}_{k'} \exp \left( \frac{\Delta V_{i k'}}{g_0 I_{sp}^s} \right) \qquad \forall k \in D_s, i \in T \label{eq:linear const6}\\
    &\left[\sum_{p \in D_s^{k'}} \left( u_p - m_{k'}^{\text{dry},s} \sum_{j \in T} z_{pj}^p \right) + m_{k'}^{\text{dry},s} + m_{k'}^{\text{dry},d} \right] \phi_{k'} \leq m_{\text{max}}^{l} \qquad \forall k' \in D_0 \label{eq:linear const7}
 \end{align}

As the problem becomes standard MILP, this can be solved by conventional MILP solving techniques such as the branch-and-bound method \cite{gurobi2025}.

\subsection{Nonlinear Programming with Fixed Routing Solution}
Although the routing solution is fixed, the global optimum cannot be determined analytically because of the nonlinearity of $\phi$ and $\Delta V$ functions. Many literature proposed various numerical methods for the NLP such as quasi-Newton and interior-point methods. The L-BFGS-B algorithm is ideal for this study because it efficiently handles bound-constrained optimization problems by leveraging the analytic first-order derivative to approximate the Hessian, avoiding the need for complex second-derivative calculations \cite{Zhu_1997_L-BFGS-B}. Its ability to efficiently manage highly nonlinear problems with simple bounds and its robust convergence properties make it both practical and effective for achieving accurate solutions.

Consider a route departing from the depot, visiting $n$ satellites, delivering each payload, and returning to the depot. Let the dry mass of servicer be $m^{\text{dry},s}$, payload mass for each satellite be $m^L_i$ ($i=\{1, 2, ..., n\}$), required $\Delta V$ to travel from previous node to next node $i$ be $\Delta V_i$ ($i=\{1, 2, ..., n + 1\}$). Then, the required propellant mass ($m^p$) becomes as follows.

\begin{align}
    m^p(\vb*{\eta}) &= \left(\left(\left( m^{\text{dry},s} \exp\left[\frac{\Delta V_{n+1}}{g_0 I_{sp}^s}\right] + m_n^L \right)\exp\left[\frac{\Delta V_{n}}{g_0 I_{sp}^s}\right] + m_{n-1}^L \right) \exp\left[\frac{\Delta V_{n-1}}{g_0 I_{sp}^s}\right] + \cdots  + m_1^L\right)\exp\left[\frac{\Delta V_{1}}{g_0 I_{sp}^s}\right] - m^{\text{dry},s} - \sum_{i = 1}^n m^L_i\\
    &= m^{\text{dry},s} \exp\left[\frac{\sum_{i=1}^{n+1}\Delta V_{i}}{g_0 I_{sp}^s}\right] + \sum_{i=1}^n m_i^L \exp\left[ \frac{\sum_{j=1}^i\Delta V_j}{g_0 I_{sp}^s}\right] - m^{\text{dry},s} - \sum_{i=1}^n m^L_i \label{eq:propellant mass form 1}
\end{align}

Note that only $\Delta V_1$ and $\Delta V_{n+1}$ vary when the depot's location is changed. In other words, all other $\Delta V_i$ for $i = \{2, ..., n\}$ are constant. So, we can rewrite the Eq.~\eqref{eq:propellant mass form 1} as follows.

\begin{align}
    m^p(\vb*{\eta}) &= A \exp \left[ \frac{\Delta V_{1}(\vb*{\eta}) + \Delta V_{n+1}(\vb*{\eta})}{g_0 I_{sp}^s} \right] + B \exp \left[ \frac{\Delta V_{1}(\vb*{\eta})}{g_0 I_{sp}^s} \right] - C
    \label{eq:propellant mass with ABC}
\end{align}
where $A = m^{\text{dry,s}} \exp\left[\frac{\sum_{i=2}^{n}\Delta V_{i}}{g_0 I_{sp}^s}\right]$, $B = \sum_{i=1}^n m_i^L \exp \left[\frac{\sum_{j=2}^i \Delta V_j}{g_0 I_{sp}^s}\right]$, and $C = m^{\text{dry,s}} + \sum_{i=1}^n m_i^L$, which are all constants.

Let $m_{r,k}^p$ be the required propellant for depot $k$'s route $r$ ($\in R_{k}$) where $R_{k}$ is the set of used routes at depot $k$ given by the routing solution of $P_L$. Then, the NLP to minimize the total amount of propellant in terms of the EMLEO mass is given as follows:

($P_{NL}$) Nonlinear programming with fixed routing solution
\begin{align}
    \min_{\vb*{\eta}} J(\vb*{\eta}) = \min_{\vb*{\eta}_{k'}} \sum_{k' \in D_0} \sum_{r \in R_{k'}} m_{r,k}^p(\vb*{\eta}_{k'}) \phi_{k'}(\vb*{\eta}_{k'})
\end{align}
subject to
\begin{align}
    r_{\text{min}} &\leq a_{k'} \qquad \forall k' \in D_0
\end{align}
where $r_{\text{min}}$ denotes the minimum radius (semimajor axis) that the depots can have.

The partial derivative of $J$ with respect to $\vb*{\eta}_{k'}$ is given as:
\begin{align}
    \frac{\partial J}{\partial \vb*{\eta}_{k'}} = \sum_{r \in R_{k'}} \left( \frac{\partial m_{r,k}^p}{\partial \vb*{\eta}_{k'}}\phi_{k'} + m_{r,k}^p \frac{\partial\phi_{k'}}{\partial \vb*{\eta}_{k'}} \right) \qquad \forall k' \in D_0
\end{align}

For notational convenience, the subscript $k'$ is omitted without loss of generality because the derivative that will be derived below can be applied to any $k'$, and they are independent of each other. Mathematically,

\begin{align}
    \pdv{m^p_{r,k_1}}{\vb*{\eta}_{k_2}} = \pdv{\phi_{k_1}}{\vb*{\eta}_{k_2}} = 0 \qquad \text{if } k_1 \neq k_2
\end{align}

Also, we use the canonical unit, i.e., $\mu = 1~\text{DU}^3/\text{TU}^2$ and $1~\text{DU}$ is the orbital radius of the constellation (e.g., 26,560 km for GPS constellation). Then, the derivative of $\phi(a)$ with respect to $a$ is as follows:

\begin{align}
    \pdv{\phi}{a} &=  \exp \left[ 
    \frac{ \sqrt{\frac{1}{a}} - \sqrt{\frac{2}{a} - \frac{2}{a + r_0}}}{g_0 \ I_{sp}^d} + \frac{\sqrt{\frac{2}{r_0} - \frac{2 }{a + r_0}} -\sqrt{\frac{1}{r_0}}}{g_0 I_{sp}^l} \right] \left( \frac{-\frac{1}{2 a^{3/2}} - \frac{-\frac{1}{a^2} + \frac{1}{(a + r_0)^2}}{\sqrt{\frac{2}{a} - \frac{2}{a + r_0}}}}{g_0 I_{sp}^d} + \frac{1}{g_0 I_{sp}^l (a+r_0)^2 \sqrt{\frac{2}{r_0} - \frac{2}{a+r_0}}} \right)
\end{align}
while the derivatives with respect to the inclination and RAAN are all 0. 

\begin{align}
    \pdv{\phi}{i}=\pdv{\phi}{\Omega} = 0
\end{align}

Now, consider the derivative of $m^p_{r,k}$ with respect to $a_{k}$. Similarly, for notational convenience, the subscripts $r$ and $k$ are omitted without loss of generality because the derivative that will be derived below can be applied to any route $r$ of depot $k$. From Eq.~\eqref{Eq:Edelbaum_dV_min}, $\Delta V_1$ and $\Delta V_{n+1}$ are the function of depot's orbital elements, $\vb*{\eta}$, given as follows.

\begin{align}
    \Delta V_1(\vb*{\eta}) &= \sqrt{\frac{1}{a} + \frac{1}{a_1} - \frac{2}{\sqrt{a a_1}} \cos \left[\frac{\pi}{2} \min \qty{\Delta \theta_1(\vb*{\eta}), 2} \right]}\label{eq:DV_1}\\
    \Delta V_{n+1}(\vb*{\eta}) &= \sqrt{\frac{1}{a} + \frac{1}{a_{n+1}} - \frac{2}{\sqrt{a a_{n+1}}} \cos \left[\frac{\pi}{2} \min\qty{\Delta \theta_{n+1}(\vb*{\eta}), 2} \right]}
    \label{eq:DV_n+1}
\end{align}
where $a_1$ and $a_{n+1}$ are the semimajor axis of the first- and last-visiting satellite, while $\Delta \theta_1$ and $\Delta \theta_{n+1}$ are the angles between the depot's orbital plane and the first and last satellite's orbital plane, respectively. Thus, the derivative of $m^p$ with respect to $\vb*{\eta}$ is given as follows:

\begin{align}
    \pdv{m^p}{\vb*{\eta}} = A \exp \qty[ \frac{\Delta V_{1} + \Delta V_{n+1}}{g_0 I_{sp}^s} ] \frac{1}{g_0 I_{sp}^s} \qty(\pdv{\Delta V_1}{\vb*{\eta}} + \pdv{\Delta V_{n+1}}{\vb*{\eta}}) + B \exp \qty[ \frac{\Delta V_{1}}{g_0 I_{sp}^s}] \frac{1}{g_0 I_{sp}^s} \pdv{\Delta V_1}{\vb*{\eta}}
\end{align}
where 

\begin{align}
    \pdv{\Delta V_1}{a} &= \frac{1}{2 \Delta V_1} \qty(-\frac{1}{a^2} + \frac{1}{\sqrt{a^3 a_1}} \cos \qty[\frac{\pi}{2} \min \qty{\Delta \theta_1, 2}])\\
    \pdv{\Delta V_1}{i} &= \frac{1}{\Delta V_1 \sqrt{a a_1}} \sin \qty[\frac{\pi}{2} \min \qty{\Delta \theta_1, 2}] \frac{\pi}{2} \pdv{\Delta \theta_1}{i}\\
    \pdv{\Delta V_1}{\Omega} &= \frac{1}{\Delta V_1 \sqrt{a a_1}} \sin \qty[\frac{\pi}{2} \min \qty{\Delta \theta_1, 2}] \frac{\pi}{2} \pdv{\Delta \theta_1}{\Omega}\\
    \pdv{\Delta V_{n+1}}{a} &= \frac{1}{2 \Delta V_{n+1}} \qty(-\frac{1}{a^2} + \frac{1}{\sqrt{a^3 a_{n+1}}} \cos \qty[\frac{\pi}{2} \min \qty{\Delta \theta_{n+1}, 2}])\\
    \pdv{\Delta V_{n+1}}{i} &= \frac{1}{\Delta V_{n+1} \sqrt{a a_{n+1}}} \sin \qty[\frac{\pi}{2} \min \qty{\Delta \theta_{n+1}, 2}] \frac{\pi}{2} \pdv{\Delta \theta_{n+1}}{i}\\
    \pdv{\Delta V_{n+1}}{\Omega} &= \frac{1}{\Delta V_{n+1} \sqrt{a a_{n+1}}} \sin \qty[\frac{\pi}{2} \min \qty{\Delta \theta_{n+1}, 2}] \frac{\pi}{2} \pdv{\Delta \theta_{n+1}}{\Omega}
\end{align}
\begin{align}
    \pdv{\Delta \theta_1}{i} &= \begin{cases}
        -\frac{1}{\sqrt{1 - \cos^2 \Delta \theta_1}} \qty(\cos i \sin i_1 \cos \qty(\Omega - \Omega_1) - \sin i \cos i_1), & \Delta \theta_1 \leq 2\\
        0, & \text{otherwise}
    \end{cases}\\
    \pdv{\Delta \theta_1}{\Omega} &= \begin{cases}
        \frac{1}{\sqrt{1 - \cos^2 \Delta \theta_1}} \sin i \sin i_1 \sin \qty(\Omega - \Omega_1), & \Delta \theta_1 \leq 2\\
        0, & \text{otherwise}
    \end{cases}\\
    \pdv{\Delta \theta_{n+1}}{i} &= \begin{cases}
        -\frac{1}{\sqrt{1 - \cos^2 \Delta \theta_{n+1}}} \qty(\cos i \sin i_{n+1} \cos \qty(\Omega - \Omega_{n+1}) - \sin i \cos i_{n+1}), & \Delta \theta_{n+1} \leq 2\\
        0, & \text{otherwise}
    \end{cases}\\
    \pdv{\Delta \theta_{n+1}}{\Omega} &= \begin{cases}
        \frac{1}{\sqrt{1 - \cos^2 \Delta \theta_{n+1}}} \sin i \sin i_{n+1} \sin \qty(\Omega - \Omega_{n+1}), & \Delta \theta_{n+1} \leq 2\\
        0, & \text{otherwise}
    \end{cases}
\end{align}

Note that $i_1$/$i_{n+1}$ and $\Omega_1$/$\Omega_{n+1}$ represent the inclination and RAAN of the first and last satellite orbits, respectively. Therefore, the first derivatives of the objective function can be calculated analytically, ensuring fast computation and accuracy for the L-BFGS-B algorithm.

The computed gradients are directly passed to the L-BFGS-B optimizer, which uses them to approximate the inverse Hessian and update the decision variables (i.e., the orbital elements $\vb*{\eta}_{k'}$ of each depot) within the given bounds. In each iteration, the optimizer evaluates both the objective function $J(\vb*{\eta})$ and its gradient $\partial J/\partial \vb*{\eta}$ to decide on the search direction and step size, based on a limited-memory quasi-Newton method. Since the gradients are available analytically in our problem, we avoid numerical differentiation, which helps reduce the number of function evaluations and prevents issues like truncation error or subtractive cancellation.

Algorithm \ref{Alg:CLRPOD_solution_framework} presents the pseudocode for the solution methodology of the CLRPOD. Starting with an initial guess, the framework iteratively updates the depot locations until convergence. Notably, solving subproblem $P_L$ is accelerated using a warm start, the routing solution from the previous iteration. This approach leverages the observation that the routing solution remains unchanged if the new depot locations do not significantly differ from the previous ones.

This framework is a generalized scheme that is not limited to specific problems or algorithms, leaving plenty of room for expansion. For instance, the initial guess can be provided by a mission planner or automatically generated using various methods, such as k-means clustering or machine learning techniques. Furthermore, not only can exact solvers (e.g., Gurobi) be directly implemented to solve $P_L$, but heuristic algorithms can also be used. Similarly, any nonlinear programming method can be utilized to solve $P_{NL}$. Since the solution quality and computational effort highly depend on the choice of the initial guess and the solvers or algorithms for $P_L$ and $P_{NL}$, the mission planner should carefully select these components based on the problem and its specific requirements.

\begin{algorithm}
\caption{Solution Framework for the CLRPOD}
\label{Alg:CLRPOD_solution_framework}
\begin{algorithmic}[1]
    \State \textbf{Initialize} depot locations: $\vb*{\eta} \gets \vb*{\eta}_0$
    \State \textbf{Initialize} routing solution for warm start: $\vb*{z} \gets \texttt{None}$
    \For{$\texttt{iter} = 1$ to $\texttt{iter}_{\max}$}
        \State Solve the routing subproblem $P_L$ to obtain an updated routing solution $\vb*{z}_{\text{new}}$ using the current depot locations $\vb*{\eta}$ and the warm start $\vb*{z}$
        \State Solve the depot location subproblem $P_{NL}$ to obtain new depot locations $\vb*{\eta}_{\text{new}}$ using the updated routing solution $\vb*{z}_{\text{new}}$
        
        \If{$\|\vb*{\eta}_{\text{new}} - \vb*{\eta}\| < \texttt{tolerance}$}
            \State \textbf{Return} depot locations $\vb*{\eta}_{\text{new}}$ and routing solution $\vb*{z}_{\text{new}}$
        \EndIf
        
        \State \textbf{Update} depot locations: $\vb*{\eta} \gets \vb*{\eta}_{\text{new}}$
        \State \textbf{Update} warm start solution: $\vb*{z} \gets \vb*{z}_{\text{new}}$
    \EndFor
    \State \textbf{Return} final depot locations $\vb*{\eta}$, routing solution $\vb*{z}$, and a message: \texttt{Max iterations reached}
\end{algorithmic}
\end{algorithm}

\section{Case Study}\label{Sec:Case_study}
This section demonstrates the mathematical formulation and solution framework of the CLRPOD through a case study focusing on a servicing scenario for a satellite constellation. For simplicity, all satellites' and depots' orbits are approximated as circular orbits so that only three orbital elements (i.e., semimajor axis, inclination, and RAAN) are sufficient to determine each orbit. 
Since we assume a static problem with regular satellite servicing (i.e., no emergency or time-critical tasks), all services can be fully planned in advance. Accordingly, each depot carries the necessary amount of fuel to its optimized location, as determined by the proposed solution framework.
A servicer corresponding to each depot departs from the depot, visits multiple satellites, performs in-space servicing, and returns to the depot. If it is necessary, the servicer will depart again for other satellites. 

To model realistic conditions, several constraints are imposed. First, the maximum launch mass is limited, and the minimum radius of the launch vehicle's parking orbit is fixed. Furthermore, the number of routes (departures) for each servicer is capped, as frequent depot revisits increase system complexity and raise the risk of mission failure. While the proposed framework allows for different dry masses for depots and servicers, this case study assumes uniform masses for simplicity, as the goal is to demonstrate the methodology rather than implement specific configurations.

Table 1 presents the parameter values used in this case study. The launch vehicle is assumed to be the Ariane 64 \cite{Ariane6_2021}. The specific impulse of the servicer is derived from \cite{Sarton_on-orbit_2022}, while the other parameters are based on \cite{shimane_orbital_2024}. For convenience, the radius of the virtual orbit for EMLEO and the minimum orbit radius for the depot are assumed to be the same ($r_0 = r_{\text{min}}$). Both the launch vehicle and the depot utilize a chemical propulsion system, whereas the servicer employs an electric propulsion system, which has a higher specific impulse. Additionally, for medium to large-scale problems, the optimal solution for the $P_L$ may not be found within the preset computation time limit, which is set at 100 seconds for each MILP in this study. Therefore, the best feasible solution obtained before reaching this time limit is used as the routing solution for the $P_L$.

The problem and solution framework are implemented in Python using the GUROBI optimizer \cite{gurobi2025}, executed on a 13th generation Intel i9 processor with a 3.0 GHz clock speed and 80 GB of RAM, operating on Windows 11.

\begin{table}[hbt!]
    \centering
    \caption{\label{tab:table_param} Parameters for GPS case study and numerical experiments}
    \begin{tabular}{lr}
    \hline
    Parameter                                                                       & Value \\
    \hline
    Maximum launch vehicle mass in the EMLEO perspective, $m_{\text{max}}^l$, kg    & 12,950\\
    Virtual orbit radius for EMLEO calculation, $r_0$, km                           & 7,000\\
    Minimum orbit radius available for the depot, $r_{\text{min}}$, km              & 7,000\\
    Standard acceleration of gravity, $g_0$, m/s\textsuperscript{2}                 & 9.81\\
    Launch vehicle's specific impulse, $I_{sp}^l$, s                                & 457\\
    Depot's specific impulse, $I_{sp}^d$, s                                         & 320\\
    Servicer's specific impulse, $I_{sp}^s$, s                                      & 1,790\\
    Maximum number of routes per depot, $n_v$, -                                           & 2\\
    Payload mass of satellite $j$, $m_j^L$, kg                                      & 100\\
    Servicer dry mass of depot $k$, $m_k^{\text{dry},s}$, kg                        & 500\\
    Depot $k$'s dry mass, $m_k^{\text{dry},d}$, kg                                  & 1500\\
    \hline
    \end{tabular}
    \label{tab:my_label}
\end{table}

\subsection{Test Case with GPS Constellation}
This case study focuses on the servicing of the GPS satellites, consisting of 18 satellites ($n_t = 18$) and three depots ($n_d = 3$). The GPS satellites are positioned in MEO at a radius of approximately 26,560 km. They have nearly circular orbits, with similar inclinations of around $55^\circ$, and are separated by $60^\circ$ in RAAN. Detailed orbital elements are provided in Appendix A (Table \ref{tab:table_sat}). Note that eccentricity and the argument of perigee are not included, as we assume the orbits are circular. Additionally, the actual RAAN of depots may vary due to J2 perturbation; see Appendix B for the corresponding sensitivity analysis.
The initial depot locations are determined using k-means clustering, where the centroids of the resulting clusters serve as the initial guesses for the depot positions.

\begin{figure}[hbt!]
    \centering
    \subfloat[]{
        \includegraphics[width=0.45\textwidth]{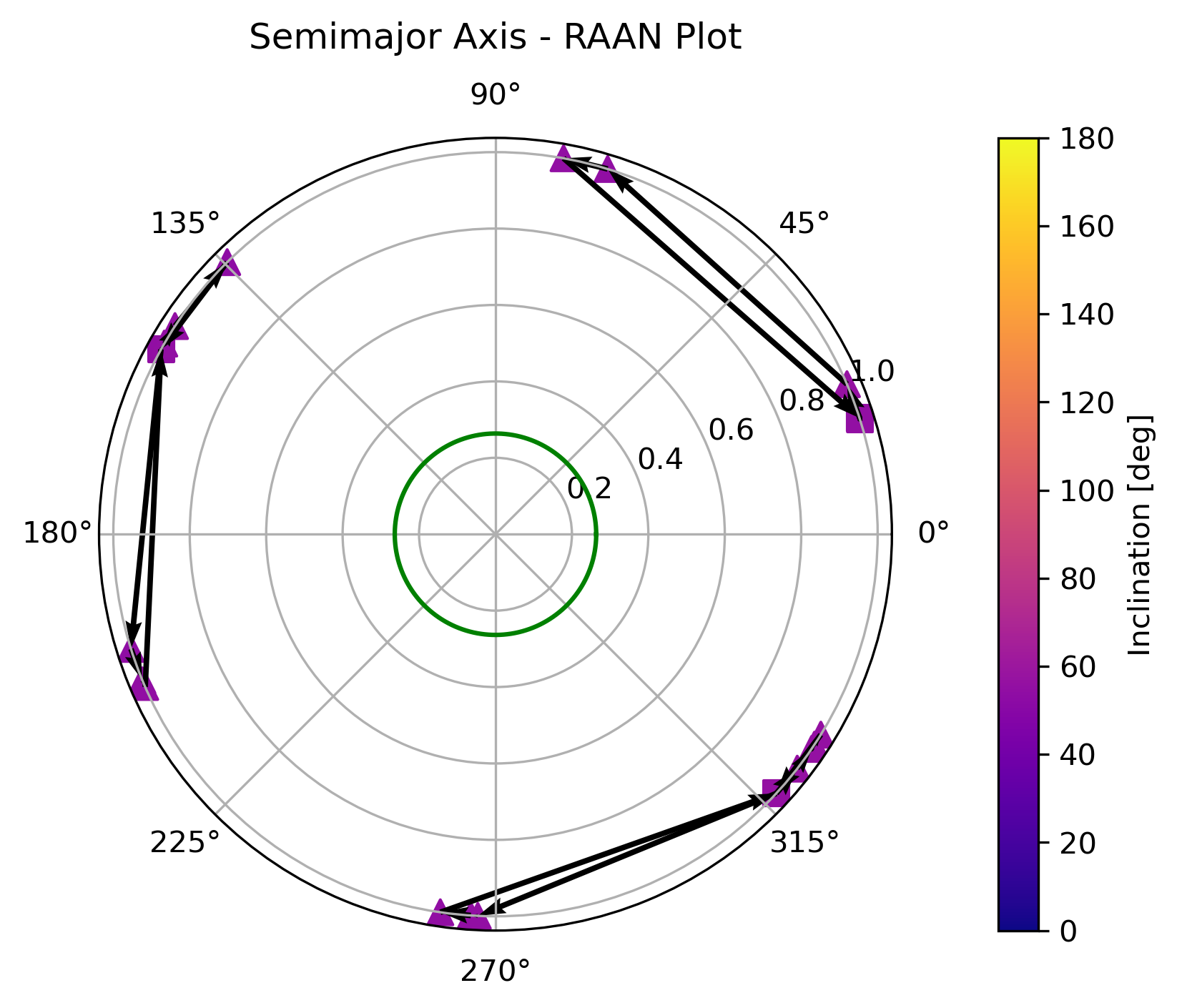}
        \label{fig:a_raan_gps1}
    }
    \hfill
    \subfloat[]{
        \includegraphics[width=0.45\textwidth]{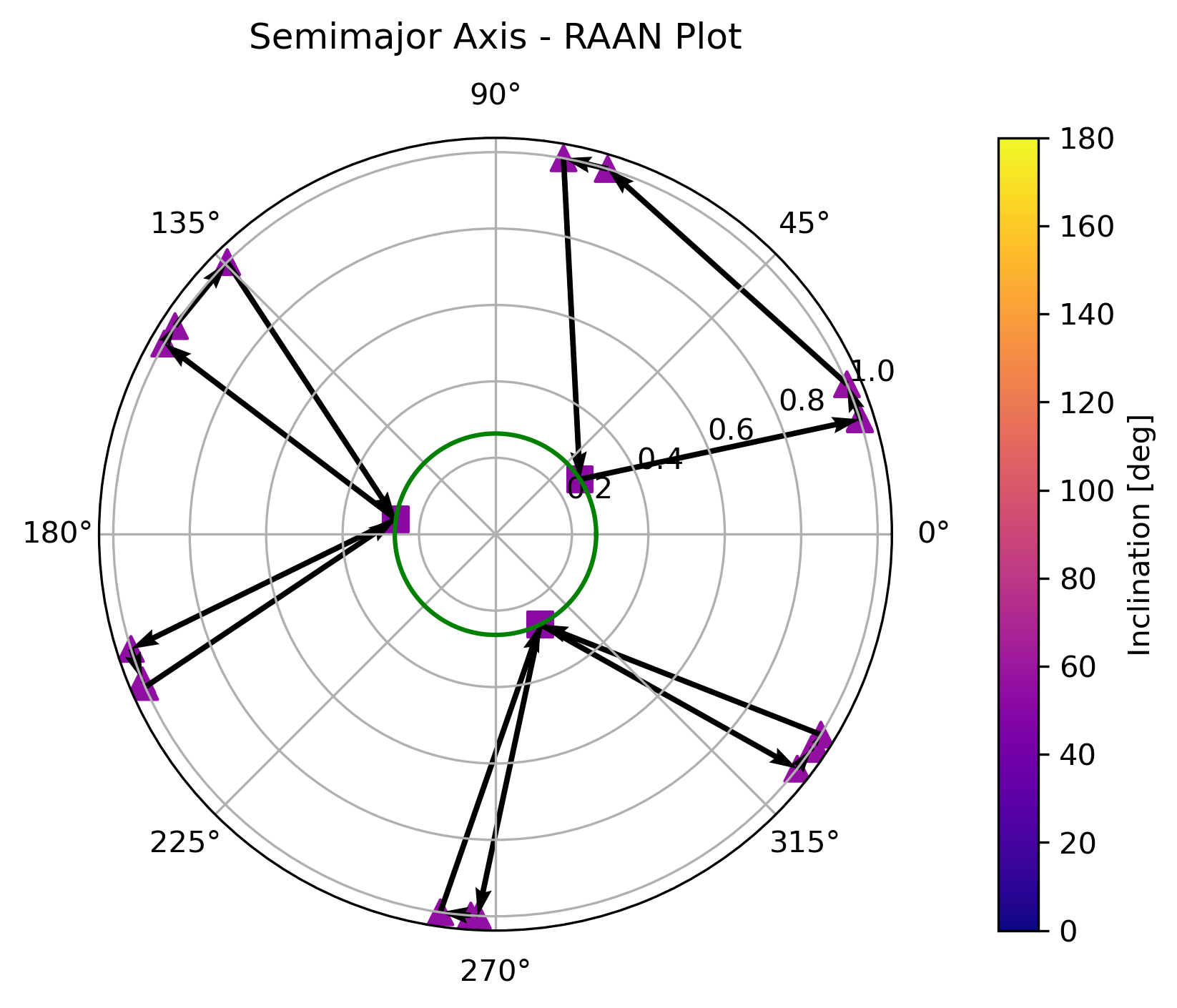}
        \label{fig:a_raan_gps2}
    }
    \caption{Comparison of the in-space servicing solutions for this case study ($a$-$\Omega$ plot): a) Solution with initial depot location guesses; b) Solution with final depot locations from the proposed framework.}
    \label{fig:a_raan_gps}
\end{figure}

\begin{figure}[hbt!]
    \centering
    \subfloat[]{
        \includegraphics[width=0.45\textwidth]{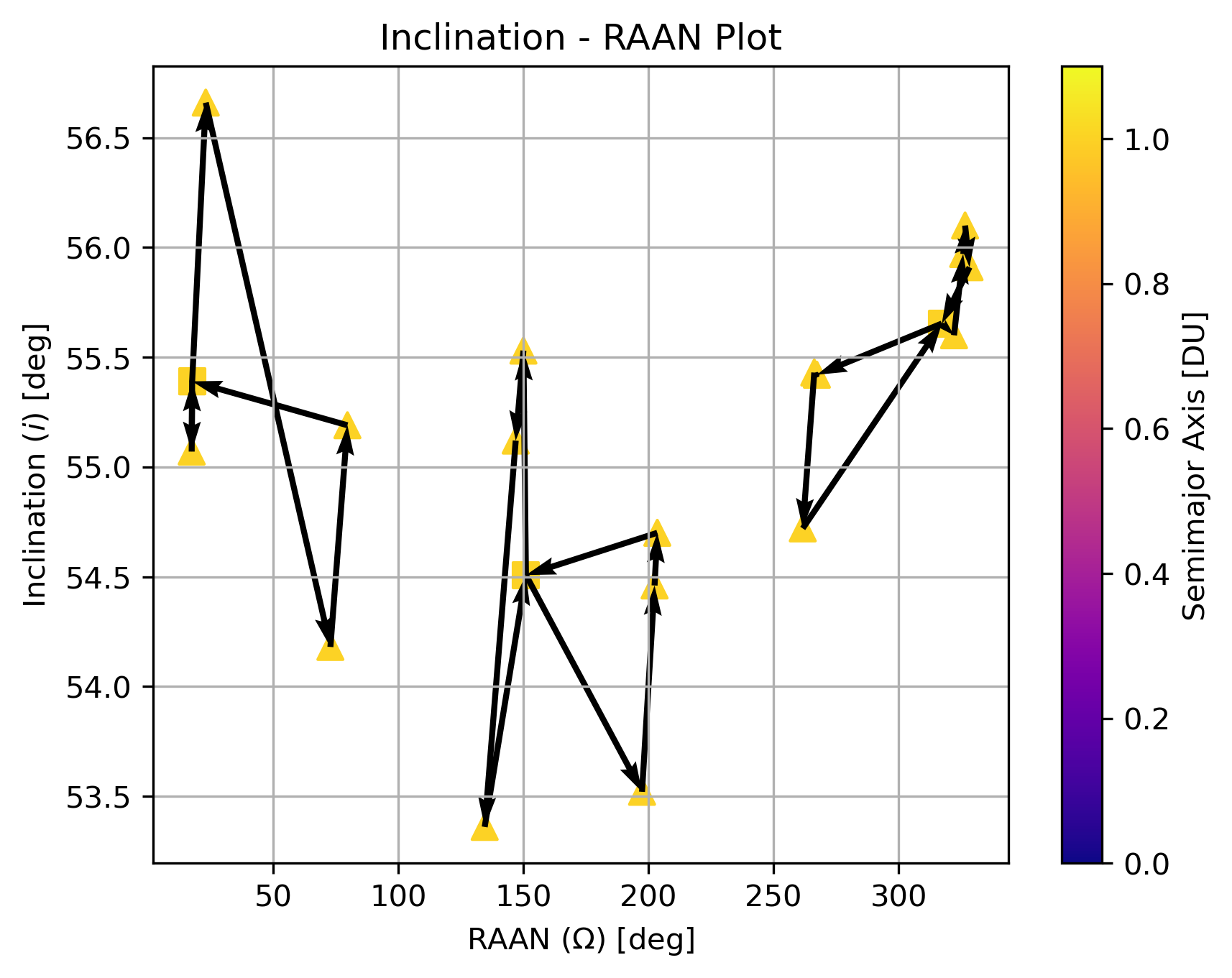}
        \label{fig:incl_raan_gps1}
    }
    \hfill
    \subfloat[]{
        \includegraphics[width=0.45\textwidth]{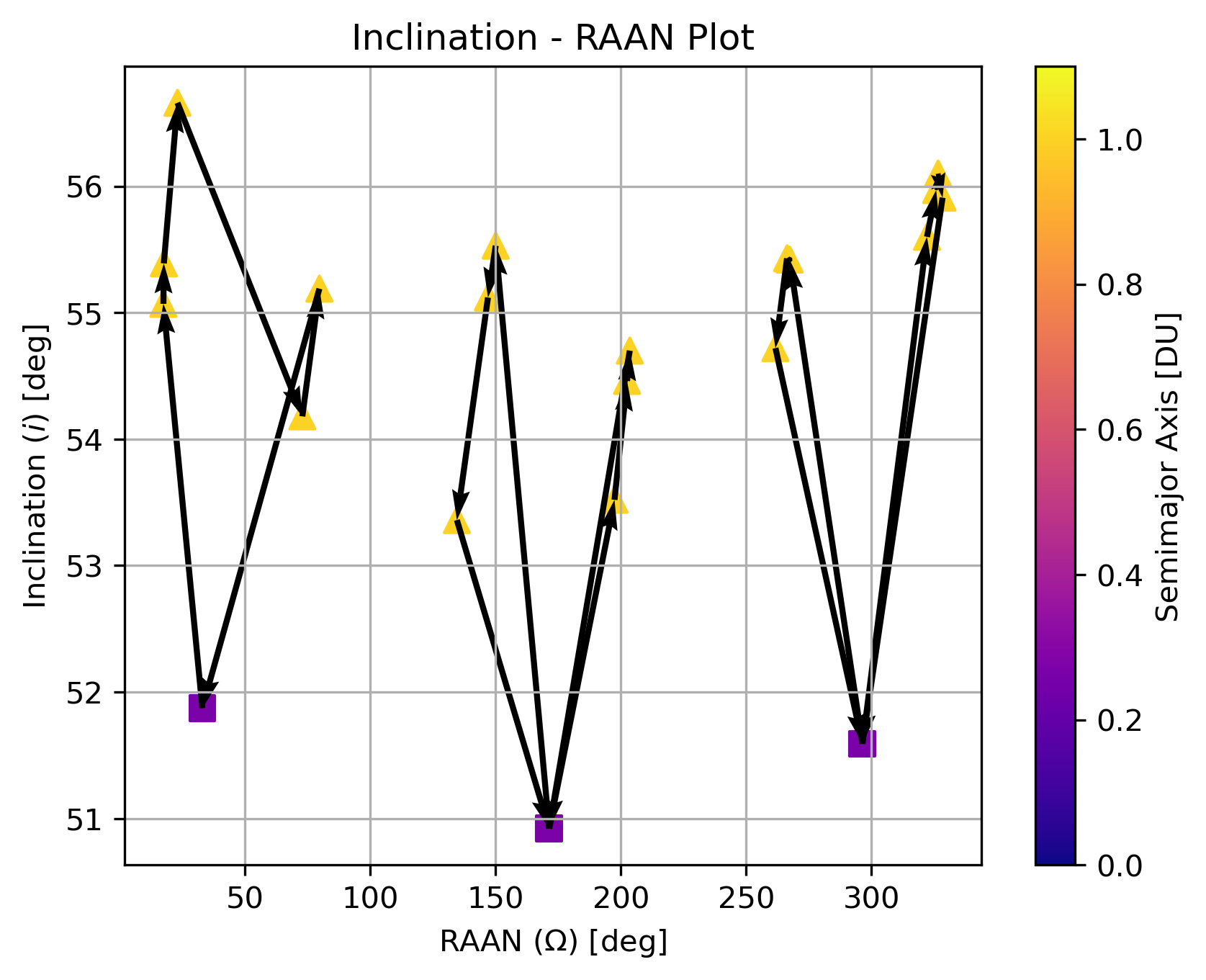}
        \label{fig:incl_raan_gps2}
    }
    \caption{Comparison of the in-space servicing solutions for this case study ($i$-$\Omega$ plot): a) Solution with initial depot location guesses; b) Solution with final depot locations from the proposed framework.}
    \label{fig:incl_raan_gps}
\end{figure}

Figures \ref{fig:a_raan_gps} and \ref{fig:incl_raan_gps} illustrate the comparisons between the solution using the initial depot locations and the solution derived from the CLRPOD framework. In these figures, squares represent depots, triangles denote the GPS satellites, black arrows indicate the visiting sequence for servicing, and the green circle in Fig. \ref{fig:a_raan_gps} marks the minimum radius allowed for the depots. 
Note that the initial depot locations obtained from k-means clustering may not appear well-centered with respect to the satellites they ultimately serve. This is due to the fundamental difference between k-means clustering and routing optimization objectives. While k-means minimizes intra-cluster variance, the routing optimization seeks to minimize overall propellant consumption. As a result, satellites may be serviced by depots outside their original k-means clusters. Nevertheless, in certain configurations, as seen in the numerical results, the optimized depot locations can coincide with the initial k-means centroids, resulting in negligible or zero EMLEO reductions.

Table \ref{tab:case_study_result} summarizes the case study results. Because the specific impulses of the launch vehicle and depot are significantly lower than that of the servicer, the depots are located at the lowest possible altitudes to minimize overall propellant consumption. 
However, the placement of depots at the lowest altitude is not universally optimal. Higher-altitude depot locations may yield better overall efficiency depending on mission-specific factors such as satellite distribution, servicing sequence, and propulsion characteristics.
Additionally, the depots are positioned with inclinations of approximately $50^\circ$ instead of $55^\circ$. This adjustment is necessary because the angle difference between orbital planes ($\Delta \theta$) is a function of both inclination and RAAN. By fine-tuning both the inclinations and RAAN, we can reduce the plane-change angle, directly affecting the amount of propellant consumed. As a result, the total propellant required in EMLEO is reduced from 7,773.982 kg to 4,906.056 kg, yielding a 36.891\% reduction compared to the initial guess. This reduction from the initial guess quantitatively demonstrates the effectiveness of optimizing the depot locations. The solution framework converges in two iterations, resulting in a computation time of 305.923 seconds.

\begin{table}[ht]
\caption{\label{tab:case_study_result} Results summary of the GPS case study}
\centering
\begin{tabular}{lcccccccc}
\hline
\multirow{2}{*}{Depot Index} & \multicolumn{3}{c}{Initial Guess} & \multicolumn{3}{c}{Final Solution} & \multirow{2}{*}{\begin{tabular}[c]{@{}c@{}}Computation \\ Time, s\end{tabular}} \\
\cmidrule(lr){2-4} \cmidrule(lr){5-7}
                             & $a$, km & $i$, deg &$\Omega$, deg & $a$, km & $i$, deg & $\Omega$, deg &                                       &                                      \\
\hline
1                            & 26,560.32  &  55.65 & 317.28         & 7,000.00 & 51.59 & 296.41 & \multirow{3}{*}{305.923} \\
2                            & 26,572.91  &  55.39 & 17.68       & 7,000.00 & 51.87 & 33.04 &  \\
3                            & 26,560.14  &  54.51 & 151.08       & 7,000.00 & 50.92 & 171.41 &  \\ \hline
Total Propellant EMLEO, kg & \multicolumn{3}{c}{7,773.982}            & \multicolumn{3}{c}{4,906.056 (36.891 \% $\downarrow$)}  &  \\
\hline
\end{tabular}
\end{table}

\subsection{Numerical Experiments}
Comprehensive numerical experiments are conducted to evaluate the performance of the CLRPOD and its solution procedure, including the EMLEO reduction from the initial guess, the number of interactions and the computational time. In this study, 100 problem instances are generated for each of five different numbers of depots (with \( n_d = 1, 2, 3, 4, \) and \( 5 \)) and four different numbers of satellites (with \( n_t = 5, 10, 15, \) and \( 20 \)). This results in a total of 2,000 problems being solved. The semimajor axes of all satellites are randomly generated in the range from 7,000 km (LEO) to 42,164 km (GEO). Both the inclination and the RAAN are varied within the ranges of \( [0^\circ, 180^\circ] \) and \( [0^\circ, 360^\circ] \), respectively. The initial guess for the depot locations is automatically generated using k-means clustering, which gives fast and intuitive results as we already know how many depots are needed. However, as demonstrated in the earlier GPS case study, the proposed solution framework can also handle any initial guesses obtained from other methods.

Table \ref{tab:numerical_exp_result} summarizes the results of the numerical experiments. Note that, for certain problems with a small number of depots, no feasible solutions are found within the preset time limit for computation from the initial guess based on k-means clustering. The column \textit{Solve success rate} in the table refers to the percentage of problems that are successfully solved feasibly. Generally, the solve success rate tends to decrease as the number of satellites increases, while it increases with a greater number of depots. The min, max, and mean values of the total EMLEO reduction from the initial guess, the number of iterations, and the computational time for each case are calculated only based on the solved problems. As the number of satellites increases, the ratio of total propellant EMLEO reduction from the initial guess also tends to increase. At the same time, both the number of iterations and the time taken to converge increase with the number of depots and satellites. Note that the number of iterations is small for some cases because the initial depot locations determined by k-means clustering can be near-optimal.

\begin{table}[hbt!]
\centering
\caption{Numerical Experiment Results of CLRPOD}
\label{tab:numerical_exp_result}
\begin{tabular}{lccccccccccc}
\hline
\multirow{2}{*}{$n_d$} &
  \multirow{2}{*}{$n_t$} &
  \multicolumn{3}{c}{\begin{tabular}[c]{@{}c@{}}Total propellant \\ EMLEO reduction\\ from initial guess, \%\end{tabular}} &
  \multicolumn{3}{c}{Number of iterations, -} &
  \multicolumn{3}{c}{Computation time, s} &
  \multirow{2}{*}{\begin{tabular}[c]{@{}c@{}}Solve success rate, \%\end{tabular}} \\ \cline{3-11}
                   &    & min   & max   & mean  & min & max & mean & min    & max    & mean   &     \\ \hline
\multirow{4}{*}{1} & 5  & 34.32 & 51.94 & 41.73 & 1   & 3   & 1.48 & 0.25   & 1.79   & 0.53   & 83  \\
                   & 10 & --    & --    & --    & --  & --  & --   & --     & --     & --     & 0   \\
                   & 15 & --    & --    & --    & --  & --  & --   & --     & --     & --     & 0   \\
                   & 20 & --    & --    & --    & --  & --  & --   & --     & --     & --     & 0   \\ \hline
\multirow{4}{*}{2} & 5  & 0.00  & 62.68 & 30.86 & 0   & 3   & 1.27 & 0.66   & 4.17   & 2.14   & 100 \\
                   & 10 & 0.00  & 60.74 & 44.91 & 0   & 7   & 2.28 & 100.31 & 804.44 & 329.53 & 100 \\
                   & 15 & 46.28 & 46.28 & 46.28 & 1   & 1   & 1.00 & 202.57 & 202.57 & 202.57 & 1   \\
                   & 20 & --    & --    & --    & --  & --  & --   & --     & --     & --     & 0   \\ \hline
\multirow{4}{*}{3} & 5  & 0.00  & 65.93 & 5.26  & 0   & 2   & 0.24 & 1.06   & 9.46   & 2.88   & 100 \\
                   & 10 & 0.00  & 60.66 & 35.72 & 0   & 5   & 1.93 & 100.21 & 605.02 & 294.83 & 100 \\
                   & 15 & 0.00  & 59.01 & 44.62 & 0   & 6   & 2.50 & 100.53 & 704.76 & 353.54 & 96  \\
                   & 20 & 31.55 & 31.55 & 31.55 & 3   & 3   & 3.00 & 402.10 & 402.10 & 402.10 & 1   \\ \hline
\multirow{4}{*}{4} & 5  & 0.00  & 0.00  & 0.00  & 0   & 0   & 0.00 & 1.81   & 15.61  & 6.36   & 100 \\
                   & 10 & 0.00  & 56.23 & 16.95 & 0   & 4   & 0.89 & 100.19 & 504.30 & 190.10 & 100 \\
                   & 15 & 0.00  & 59.60 & 38.30 & 0   & 6   & 2.33 & 100.30 & 702.51 & 335.56 & 100 \\
                   & 20 & 0.00  & 62.85 & 42.94 & 0   & 6   & 2.98 & 100.58 & 708.54 & 401.51 & 91  \\ \hline
\multirow{4}{*}{5} & 5  & 0.00  & 0.00  & 0.00  & 0   & 0   & 0.00 & 0.96   & 7.80   & 2.32   & 100 \\
                   & 10 & 0.00  & 55.62 & 5.59  & 0   & 5   & 0.29 & 100.19 & 601.41 & 129.84 & 100 \\
                   & 15 & 0.00  & 54.64 & 22.81 & 0   & 5   & 1.44 & 100.35 & 604.45 & 245.59 & 100 \\
                   & 20 & 0.00  & 59.69 & 34.91 & 0   & 6   & 2.38 & 100.59 & 706.73 & 340.97 & 100 \\ \hline
\end{tabular}
\end{table}

\section{Conclusions}\label{Sec:Conclusions}
This paper presents a continuous location routing problem for orbital depots to support in-space servicing missions. One uniqueness of this problem is that depot locations are chosen from a continuous space instead of from a discrete set of candidates, which was commonly assumed in the literature. The optimization problem minimizes total propellant consumption by simultaneously optimizing the orbital elements of the depots and the sequence in which the satellites are visited. Due to the nonlinear and discrete nature of this problem, directly solving the problem is challenging. Therefore, we introduce a solution framework that decouples the problem into mixed-integer linear programming and nonlinear programming to solve iteratively. A case study involving the GPS constellation, along with numerical experiments, demonstrates the performance characteristics and effectiveness of the proposed formulation and solution framework.

Future research could focus on developing methods that handle infeasible initial guesses or consider multiple propulsion systems. Exploring more practical scenarios, such as considering the available time window for each satellite or planning with stochastic failure of satellites/servicers, can also be a promising further study. Additionally, extending the application of this work to cis-lunar space or deep space missions could provide valuable insights for upcoming space campaigns.

\section*{Appendix A: Data Table Used in Case Study} 
\setcounter{table}{0}
\renewcommand{\thetable}{A\arabic{table}}
\begin{table}[hbt!]
\caption{\label{tab:table_sat} Selected orbital elements of the GPS satellites}
\centering
\begin{tabular}{lccc}
\hline
Satellite Index & $a$, km & $i$, deg & $\Omega$, deg \\
\hline
1 & 26,560.36 & 55.53 & 150.07 \\
2 & 26,560.46 & 54.18 & 72.93 \\
3 & 26,561.19 & 55.12 & 146.99 \\
4 & 26,561.01 & 55.42 & 267.35 \\
5 & 26,560.44 & 55.07 & 17.50 \\
6 & 26,560.92 & 55.91 & 328.36 \\
7 & 26,572.91 & 55.39 & 17.68 \\
8 & 26,560.09 & 55.97 & 325.81 \\
9 & 26,559.72 & 54.70 & 203.57 \\
10 & 26,560.77 & 55.43 & 266.30 \\
11 & 26,560.02 & 53.36 & 134.59 \\
12 & 26,559.80 & 56.10 & 326.58 \\
13 & 26,559.86 & 54.46 & 202.48 \\
14 & 26,559.18 & 55.19 & 79.74 \\
15 & 26,559.54 & 54.72 & 261.70 \\
16 & 26,560.12 & 56.66 & 23.12 \\
17 & 26,560.35 & 53.52 & 197.47 \\
18 & 26,560.20 & 55.60 & 322.15 \\
\hline
\end{tabular}
\end{table}

\section*{Appendix B: Sensitivity Analysis on RAAN Drift due to J2 Perturbation}
In this study, all orbits are assumed to be fixed, without considering J2 perturbation. However, in reality, the J2 effect causes the changes in the RAAN depending on orbital altitude. For example, the RAAN of a GPS orbit changes by approximately $-0.0387^\circ/\text{day}$, whereas an orbit with a radius of 7,000 km experiences a drift of approximately $-4.121^\circ/\text{day}$.

Given the use of a low-thrust propulsion system, each orbital transfer is assumed to take approximately 10--50 days, and the overall mission (servicing 3 to 4 satellites per route) would span roughly 100 days. Therefore, the impact of J2 perturbation is not negligible, and a depot's RAAN may drift significantly throughout the mission.

To evaluate the sensitivity of the optimization results to RAAN drift, we conducted a sensitivity analysis by varying the RAAN of the depot orbit. Figure~\ref{fig:raan_sensitivity} presents the variation in the objective value (total propellant consumption in terms of EMLEO) as a function of the depot’s RAAN for the GPS case study. It is important to note that we evaluate relative RAAN differences between depots and GPS satellites; thus, only the return transfer to the depot is affected by RAAN drift. The red dot in each figure represents the initial RAAN value (i.e., the case in which no RAAN drift occurs and the depot remains fixed).

The results show that the objective value changes by an average of 6.62\%, with a maximum deviation of 14.98\% from the original solution. This level of variation is considered acceptable, as the proposed method is intended for high-level strategic planning rather than detailed tactical trajectory optimization.

It is also important to note that solutions generated by the k-means clustering would be similarly affected by RAAN drift. Therefore, the observed variation does not imply that our proposed method is more sensitive to J2 effects than the k-means clustering. The performance improvements over the k-means clustering remain valid.

\setcounter{figure}{0}
\renewcommand{\thefigure}{B\arabic{figure}}
\begin{figure}[hbt!]
    \centering
    \subfloat[]{
        \includegraphics[width=0.32\textwidth]{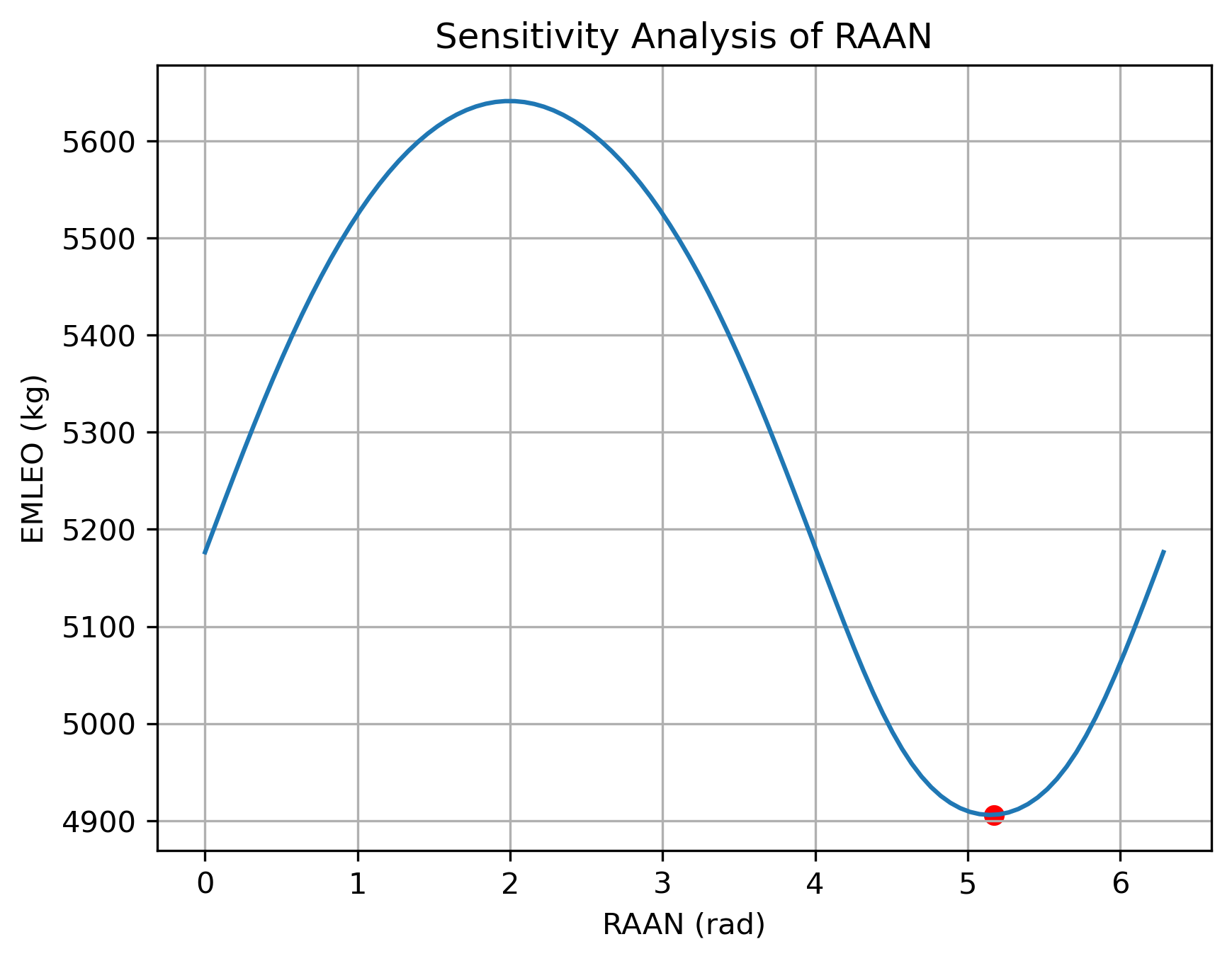}
        \label{fig:raan_sensitivity0}
    }
    \hfill
    \subfloat[]{
        \includegraphics[width=0.32\textwidth]{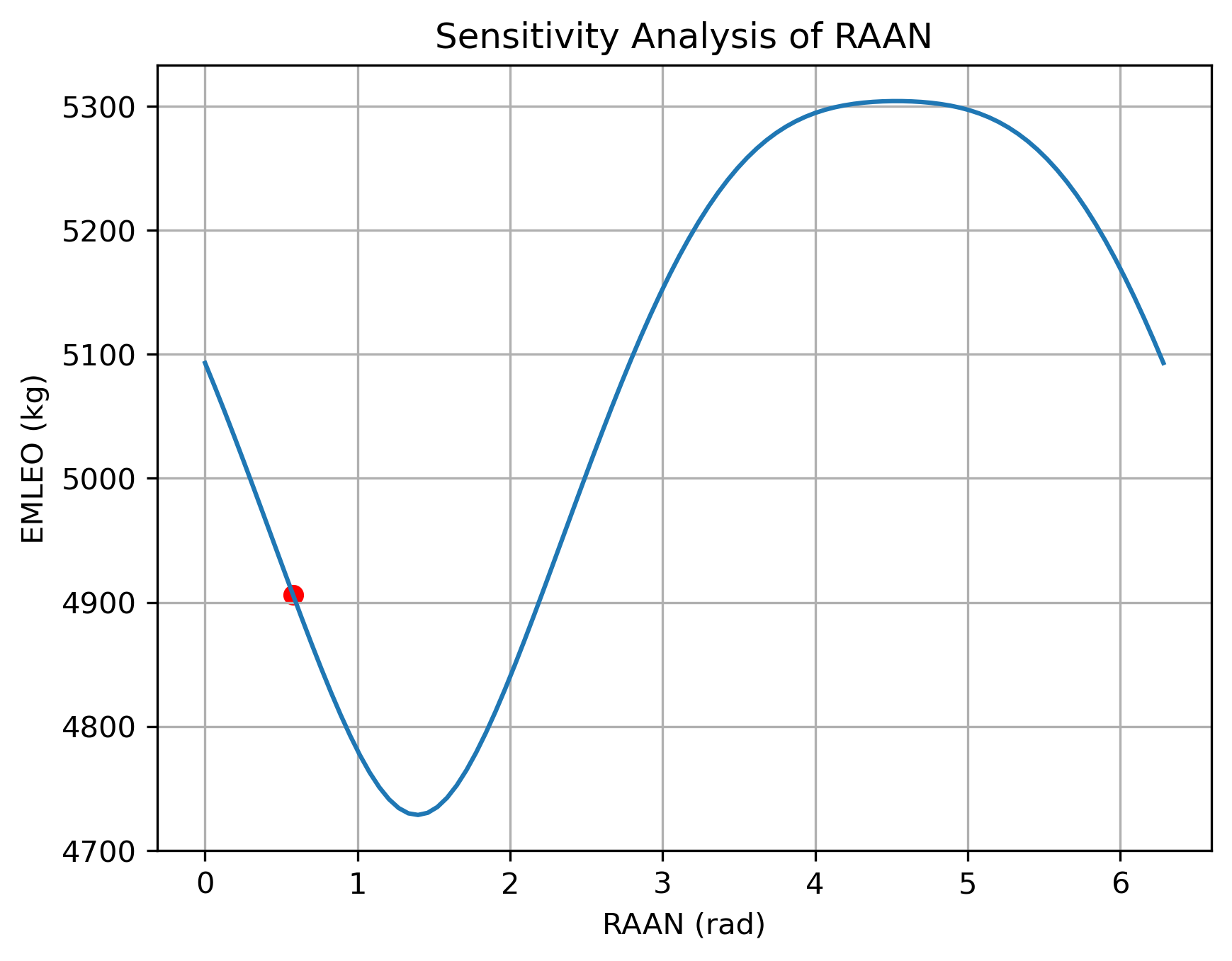}
        \label{fig:raan_sensitivity1}
    }
    \hfill
    \subfloat[]{
        \includegraphics[width=0.32\textwidth]{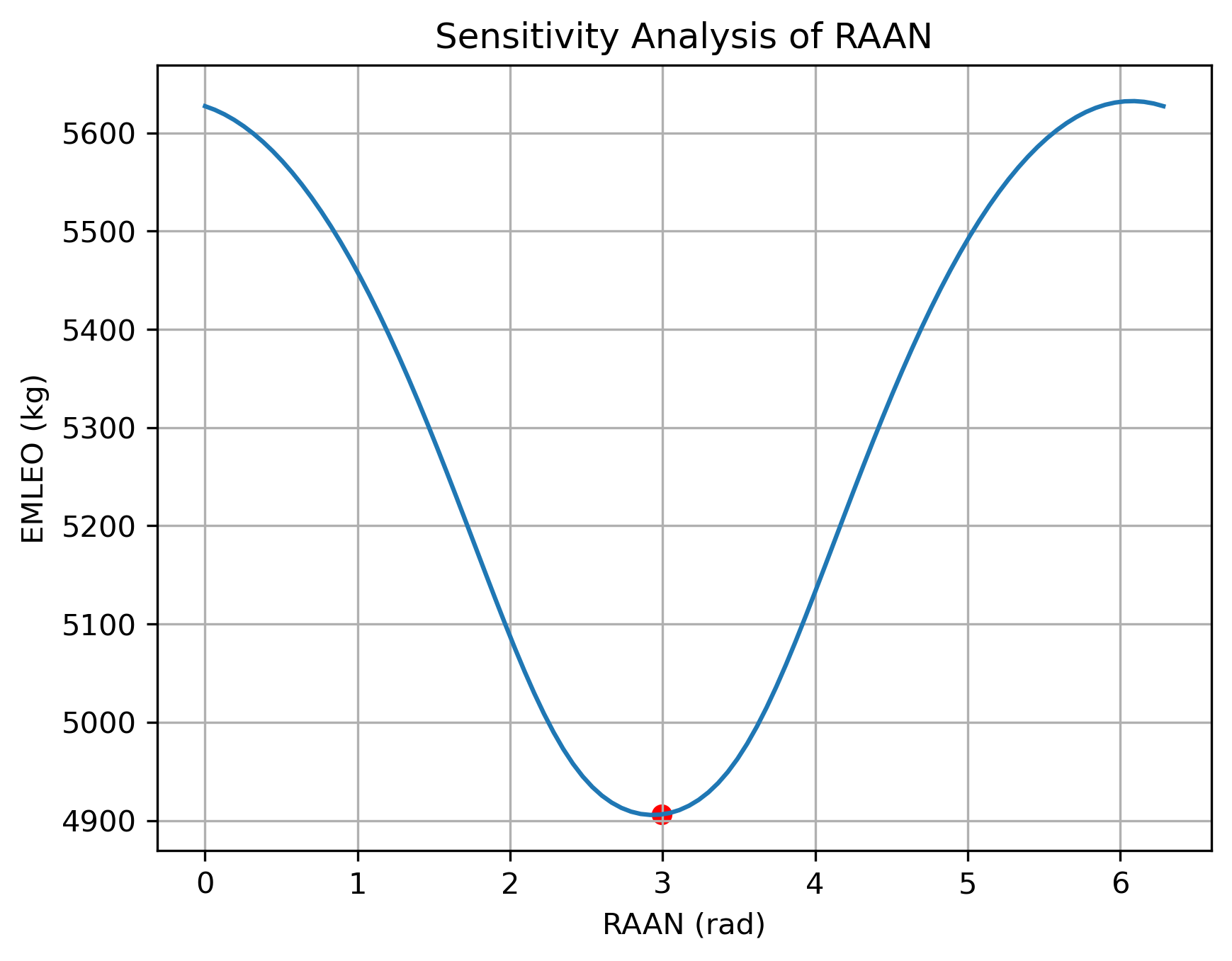}
        \label{fig:raan_sensitivity2}
    }
    \caption{EMLEO variation due to RAAN drift of: (a) Depot 1; (b) Depot 2; (c) Depot 3.}
    \label{fig:raan_sensitivity}
\end{figure}

\section*{Acknowledgments}
This work was conducted with support from the Air Force Office of Scientific Research (AFOSR), as part of the Space University Research Initiative (SURI), under award number FA9550-23-1-0723. We thank Yuri Shimane for reviewing and verifying the code. AI technologies were used for grammar checking and refinement.

\bibliography{echoi}

\end{document}